\documentclass[12pt]{amsart}

\usepackage{pstricks}\psset{xunit=2.2ex,yunit=2.2ex}
\usepackage{latexsym}
\usepackage{amssymb}

\newcommand{\excise}[1]{}

\newtheorem{thm}{Theorem}
\newtheorem{lemma}[thm]{Lemma}

\newtheorem{cor}[thm]{Corollary}
\newtheorem{prop}[thm]{Proposition}

\theoremstyle{definition}

\newtheorem{remark}[thm]{Remark}

\newtheorem{defn}[thm]{Definition}

	{\end{trivlist}}

\newenvironment{stareqn}%
	{
	 \setcounter{separated}{\value{equation}}
	 \setcounter{equation}{0}
	 \begin{eqnarray}
	}
	{\end{eqnarray}%
	 \setcounter{equation}{\value{separated}}%
	}

        {\begin{list}
                {\noindent\makebox[0mm][r]{\arabic{enumi}.}}
                {\leftmargin=5.5ex \usecounter{enumi}}
        }
        {\end{list}}

        {\begin{list}
                {\noindent\makebox[0mm][r]{(\roman{enumi})}}
                {\leftmargin=5.5ex \usecounter{enumi}}
        }
        {\end{list}}

\newenvironment{vcenteredyoung}{ \begin{array}{@{}c@{}}\\[-2.2ex]\young}
				{\end{array}}

\newcounter{separated}

\def\bem#1{\textbf{#1}}

\def\<{\langle}
\def\>{\rangle}
\def\0{{\mathbf 0}}
\def\1{{\mathbf 1}}

\def\dz{{\dot z}}

\def\GG{{\mathcal G}}
\def\fG{G}

\def\KQ{{\mathit K}\!{\mathcal Q}}

\def\NN{{\mathbb N}}

\def\PP{{\mathcal P}}

\def\ZZ{{\mathbb Z}}

\def\rr{{\mathbf r}}

\def\ww{{\mathbf w}}
\def\xx{{\mathbf x}}
\def\yy{{\mathbf y}}
\def\zz{{\mathbf z}}

\def\th{{\rm th}}

\def\length{l}

\def\K{$K$}

\def\kl{{k\!\times\!\ell}}
\def\oP{\check{P}}
\def\rp{\mathcal{RP}}

\def\rij{r_{\!ij}}

\def\oxx{\text{$\stackrel{\raisebox{-.2ex}{$\scriptscriptstyle\circ$}}%
     {\hbox{$\xx$}}\hspace{-.7ex}^{}$}}
\def\oyy{\text{$\stackrel{\raisebox{-.2ex}{$\scriptscriptstyle\circ$}}%
     {\hbox{$\yy$}}\hspace{-.7ex}^{}$}}

\def\homv{{\hspace{-.2ex}\text{\sl{Hom}}\hspace{.1ex}}}

\def\minus{\smallsetminus}

\def\goesto{\rightsquigarrow}

\def\ol#1{{\overline {#1}}}
\def\ub#1{\text{\underbar{$#1$\hspace{-.25ex}}\hspace{.25ex}}}

\def\wt#1{{\widetilde{#1}}}
\def\dem#1{{\ol \partial_{#1}}}

\def\textcross{\ \smash{\lower4pt\hbox{\rlap{\hskip4.15pt\vrule height14pt}}
                \raise2.8pt\hbox{\rlap{\hskip-3pt \vrule height.4pt depth0pt
                width14.7pt}}}\hskip12.7pt}
\def\textelbow{\ \hskip.1pt\smash{\raise2.8pt%
                \hbox{\co \hskip 4.15pt\rlap{\rlap{\char'005} \char'007}
                \lower6.8pt\rlap{\vrule height3.5pt}
                \raise3.6pt\rlap{\vrule height3.5pt}}
                \raise2.8pt\hbox{%
                  \rlap{\hskip-7.15pt \vrule height.4pt depth0pt width3.5pt}%
                  \rlap{\hskip4.05pt \vrule height.4pt depth0pt width3.5pt}}}
                \hskip8.7pt}

\font\co=lcircle10
\def\petit#1{{\scriptstyle #1}}

\def\jr{\smash{	\raise2pt\hbox{\co \rlap{\rlap{\char'005} \char'007}}
		\raise6pt\hbox{\rlap{\vrule height2pt}}
		\raise2pt\hbox{\rlap{\hskip4pt \vrule height0.4pt depth0pt
                 width2.5pt}}
		\raise2pt\hbox{\rlap{\hskip-6pt \vrule height.4pt depth0pt
                 width2.2pt}}
		\lower4pt\hbox{\rlap{\vrule height2.5pt}}}}
\def\je{\smash{\raise2pt\hbox{\co \rlap{\rlap{\char'005}
                \phantom{\char'007}}}\raise6pt\hbox{\rlap{\vrule height2pt}}
	       \raise2pt\hbox{\rlap{\hskip-6pt \vrule height.4pt depth0pt
                width2.2pt}}}}
\def\+{\smash{\lower4pt\hbox{\rlap{\vrule height12pt}}
                \raise2pt\hbox{\rlap{\hskip-6pt \vrule height.4pt depth0pt
                width12.5pt}}}}
\def\er{\smash{	\raise2pt\hbox{\co \rlap{\phantom{\rlap{\char'005}} \char'007}}
		\raise6pt\hbox{\rlap{\phantom{\vrule height2pt}}}
		\raise2pt\hbox{\rlap{\hskip4pt \vrule height0.4pt depth0pt
                 width2.5pt}}
		\raise2pt\hbox{\rlap{\phantom{%
		 \hskip-6pt \vrule height.4pt depth0pt width2.2pt}}}
		\lower4pt\hbox{\rlap{\vrule height2.5pt}}}}
\def\hor{\smash{\lower2pt\hbox{\rlap{\phantom{\vrule height10pt}}}
                \raise2pt\hbox{\rlap{\hskip-6pt \vrule height.4pt depth0pt
                width12.5pt}}}}
\def\ver{\smash{\lower2pt\hbox{\rlap{\vrule height10pt}}
                \raise2pt\hbox{\rlap{\phantom{%
		\hskip-6pt \vrule height.4pt depth0pt width12.5pt}}}}}

\def\perm#1#2{\hbox{\rlap{$\petit {#1}_{\scriptscriptstyle #2}$}}%
                \phantom{\petit 1}}
\def\phperm{\phantom{\perm w3}}

\newenvironment{tinygraph}{\begin{trivlist}\item\centering\tiny$}
			{$\end{trivlist}}

\def\normal#1{\hbox{\normalsize$#1$}}

\newenvironment{dream}[1]%
 	{
	 \def\*{\makebox[0ex]{\footnotesize$+\,$}}%
	 \begin{array}{|*{#1}{@{\ \;\:}c|}}}
	{\end{array}}

\newenvironment{pipedream}[1]%
	{
	 \begin{array}{c*{#1}{@{\ \ \;}c}}}
	{\end{array}}

\begin{document}

\mbox{}\vspace{2.48ex}
\title{Alternating formulas for \K-theoretic quiver polynomials}
\author{Ezra Miller}
\thanks{The author was partly supported by the National Science Foundation}
\address{University of Minnesota\\ Minneapolis, Minnesota}
\email{ezra@math.umn.edu}
\date{27 November 2003}

\begin{abstract}
\noindent
The main theorem here is the \K-theoretic analogue of the
cohomological `stable double component formula' for quiver functions
in \cite{quivers}.  This \K-theoretic version is still in terms of
lacing diagrams, but nonminimal diagrams contribute terms of higher
degree.  The motivating consequence is a conjecture of
Buch on the sign-alternation of the coefficients appearing in his
expansion of quiver \K-polynomials in terms of stable Grothendieck
polynomials for partitions \cite{Buch02}.
\end{abstract}

\maketitle

{}


\section*{Introduction}

\noindent

The study of combinatorial formulas for the degeneracy loci of quivers
of vector bundles with arbitrary ranks was initiated by Buch and
Fulton~\cite{BF99}.  In that paper they proved that the cohomology
classes of such degeneracy loci can be expressed as integer sums of
products of Schur polynomials evaluated on the Chern classes of the
bundles in the quiver.  After giving an explicit algorithmic (but
nonpositive) expression for the {\em quiver coefficients}\/ appearing
therein, they also conjectured a positive combinatorial formula for
them.  This conjecture was proved in \cite{quivers}, by way of three
other positive combinatorial formulas for the {\em quiver
polynomials}.

The cohomological ideas of \cite{BF99} were extended to \K-theory in
\cite{Buch02}, where the classes of the structure sheaves of the
aforementioned degeneracy loci are expressed as integer sums of
products of stable double Grothendieck polynomials for grassmannian
permutations.  Buch proved a formula for the coefficients in this
expansion \cite[Theorem~4.1]{Buch02}, and conjectured that the signs
of these coefficients alternate in a simple manner
\cite[Conjecture~4.2]{Buch02}.

The main result here is Theorem~\ref{t:formula}, which extends the
{\em stable double component formula}\/ in terms of minimal lacing
diagrams \cite[Theorem~6.20]{quivers}, as well as some combinatorial
methods surrounding it, from cohomology to \K-theory.
Theorem~\ref{t:formula} is still in terms of lacing diagrams, but
nonminimal diagrams contribute terms of higher degree.  The purpose is
to prove Buch's conjecture as a consequence (Theorem~\ref{t:alt}),
using a sign-alternation theorem of Lascoux \cite[Theorem~4]{Las01}.
The \K-theory analogue of a formula \cite[Theorem~5.5]{quivers} for
quiver polynomials in terms of the {\em pipe dreams}\/ of Fomin and
Kirillov \cite{FK96} enters along the~way~(Theorem~\ref{t:FK}).

The proof of Theorem~\ref{t:formula} generalizes a procedure suggested
by \cite{quivers} (see Remark~6.21 there), and carried out in
\cite{Yong03}, for constructing pipe dreams associated to given lacing
diagrams.  This technique is combined with those developed in
\cite{subword} for dealing with nonreduced subwords of reduced
expressions for permutations.

Buch \cite{BuchAltSign} independently arrived at the main results and
definitions here (and more) by applying general techniques of Fomin
and Kirillov \cite{FK96,FK94}.
A~special case of the sign conjecture and \K-component formula already
appeared in~\mbox{\cite{BKTY03}}.

\subsection*{Organization}
A notion of double quiver \K-polynomial is identified via a ratio
formula in Section~\ref{sec:quiver}, in analogy with the way
(cohomological) double quiver polynomials arise in \cite{quivers}.
The `pipe formula' for quiver \K-polynomials is proved in
Section~\ref{sec:pipe}, after background on nonreduced pipe dreams and
Demazure products.  The condition on nonminimal lacing diagrams that
turns out to make them occur with sign $\pm 1$ in the \K-component
formula is defined in Section~\ref{sec:lacing}.  Rank stability of
these nonminimal lacing diagrams, proved in
Section~\ref{sec:stability}, plays the same role here as it did for
the cohomological component formula in \cite{quivers}.  The stable
\K-component formula is derived in Section~\ref{sec:formula}, after
reviewing basics regarding Grothendieck polynomials and their stable
limits.  Finally, Buch's sign alternation conjecture is proved in
Section~\ref{sec:alt}.



\section{Double quiver \K-polynomials}

\label{sec:quiver}

A $\kl$ {\em partial permutation}\/ is a $\kl$ matrix~$w$ whose
entries are either $0$ or~$1$, with at most one nonzero entry in each
row or column.  Each such matrix~$w$ can be completed to a permutation
matrix---that is, with exactly one~$1$ in each row and
column---having~$w$ as its upper-left $\kl$ corner.  Viewing
permutations as lying in the union $S_\infty = \bigcup_k S_k$ of all
symmetric groups~$S_k$, there is a unique completion~$\wt w$ of~$w$
that has minimal length $\length(\wt w)$.  For any partial
permutation~$w$, we write $q = w(p)$ if the entry of~$w$ in row~$p$
and column~$q$ equals~$1$.  If~$v$ is a permutation matrix, then the
assignment $p \mapsto v(p)$ defines a permutation in~$S_\infty$.

Let $\zz = z_1,z_2,\ldots$ and $\dot\zz = \dz_1,\dz_2,\ldots$ be
alphabets.  Writing a given polynomial~$f$ in these two alphabets over
the integers~$\ZZ$ as a polynomial in $z_i$ and~$z_{i+1}$ with
coefficients that are polynomials in the other variables, the {\em
$i^\th$ Demazure operator}\/ $\dem i$ sends~$f$~to
\begin{eqnarray*}
  \dem if &=&
  \frac{z_{i+1}f(z_i,z_{i+1})-z_if(z_{i+1},z_i)}{z_{i+1}-z_i}.
\end{eqnarray*}
Let $w_0^k$ be the permutation of maximal length in~$S_k$, and write
$s_i \in S_\infty$ for the transposition switching $i$
and~\mbox{$i+1$}.  Following \cite{LS82}, the {\em double Grothendieck
polynomial}\/ for a permutation $v \in S_k$ is defined from the
``top'' double Grothendieck polynomial $\GG_{w_0^k}(\zz/\dot\zz) =
\prod_{i+j \leq k}(1-z_i/\dz_j)$ by the recursion
\begin{eqnarray*}
  \GG_{vs_i}(\zz/\dot\zz) &=& \dem i \GG_v(\zz/\dot\zz)
\end{eqnarray*}
whenever $vs_i$ is lower in Bruhat order than~$v$.  This definition is
independent of the choice of~$k$ \cite{LS82}.  If~$w$ is a partial
permutation, then set $\GG_w(\zz/\dot\zz) = \GG_{\wt w}(\zz/\dot\zz)$.

The permutation matrices~$v$ of central importance here are those
associated to the `Zelevinsky permutations' of \cite{quivers}, which
are defined as follows.  Fix a positive integer~$d$ an expression $d =
\sum_{j=0}^n r_j$ of~$d$ as a sum of $n+1$ ``ranks''~$r_j$.  Endow
each $d \times d$ permutation matrix~$v$ with a block decomposition in
which the $j^\th$ block row from the top has height~$r_j$, and the
$i^\th$ block column {\em from the right}\/ has width $r_i$.  Thus
each $d \times d$ permutation matrix~$v$ is composed of $(n+1)^2$
blocks~$B_{ji}$, each of size $r_j \times r_i$.  The matrix~$v$ is a
{\em Zelevinsky permutation}\/ as in \cite[Definition~1.7]{quivers} if
$B_{ji}$ has all zero entries whenever $i \geq j+2$, and the nonzero
entries of~$v$ proceed from northwest to southeast within every block
row or block column (so $v$ has no $1$~entry that is northeast of
another within the same block row or block column).  Pictures and
examples can be found in \cite[Section~1.2]{quivers}.

If a Zelevinsky permutation~$v$ is given, define $\rij$ to be the
number of nonzero entries of~$v$ in the union of all blocks~$B_{qp}$
for which $q \geq j$ and $p \leq i$ (that is, blocks~$B_{qp}$ weakly
southeast of~$B_{ji}$).  This results in a {\em rank array}\/ $\rr =
(\rij)_{i \geq j}$.  Since $\rr$ uniquely determines~$v$ by
\cite[Proposition~1.6]{quivers}, the notation $v = v(\rr)$ makes
sense.

Double Grothendieck polynomials for Zelevinsky permutations $v(\rr)$
are naturally written as $\GG_{v(\rr)}(\xx/\oyy)$, using two alphabets
$\zz = \xx$ and $\dot\zz = \oyy$ each of which is an ordered sequence
of $n+1$ alphabets of sizes $r_0,\ldots,r_n$ and $r_n,\ldots,r_0$,
repsectively:
\begin{eqnarray*}
  \xx = \xx^0,\ldots,\xx^n &\text{and}& \oyy = \yy^n,\ldots,\yy^0,
\\
  \makebox[0ex][r]{where\quad} \xx^j = x^j_1,\ldots,x^j_{r_{\!j}}
  &\text{and}& \yy^j = y^j_1,\ldots,y^j_{r_{\!j}}.
\end{eqnarray*}
It is convenient to think of the $\xx$ variables as labeling the rows
of the $d \times d$ grid, while the~$\yy$ variables label its columns
(see \cite[Section~2.2]{quivers} for pictures and examples).  Most
partial permutations~$w$ that occur in the sequel will have size
$r_{j-1} \times r_{\!j}$ for some~$j \in \{1,\ldots,n\}$; in that case
we consider~$\GG_w(\xx^{j-1}/\yy^j)$.

Among all $d \times d$ Zelevinsky permutations with block
decompositions determined by $d = \sum_{j=0}^n r_j$, there is a unique
one $v(\homv)$ whose rank array $\rr(\homv)$ is maximal, in the sense
that $\rij(\homv) \geq \rij$ for all other $d \times d$\/
Zelevinsky permutations~$v(\rr)$.

\begin{defn} \label{d:ratio}
The \bem{double quiver \K-polynomial} is the ratio
\begin{eqnarray*}
  \KQ_\rr(\xx/\oyy) &=& \frac{\GG_{v(\rr)}(\xx/\oyy)}
  {\GG_{v(\homv)}(\xx/\oyy)}
\end{eqnarray*}
of double Grothendieck polynomials for $v(\rr)$ and~$v(\homv)$.%
\end{defn}
\noindent

The ``ordinary'' specialization of the polynomial $\KQ_\rr(\xx/\oyy)$
appears in the \K-theoretic ratio formula \cite[Theorem~2.7]{quivers}.
It will follow from Theorem~\ref{t:FK}, below, that
$\GG_{v(\homv)}(\xx/\oyy)$ divides $\GG_{v(\rr)}(\xx/\oyy)$, so the
right hand side of Definition~\ref{d:ratio} is actually a (Laurent)
polynomial rather than simply a rational function.

\section{Nonreduced pipe dreams}

\label{sec:pipe}

A $\kl$ {\em pipe dream}\/ is a subset of the \mbox{$\kl$} grid,
identified as the set of crosses in a tiling of the $\kl$ grid by {\em
crosses}\/ $\textcross$ and {\em elbow joints}\/ $\textelbow$, as in
the following diagrams:
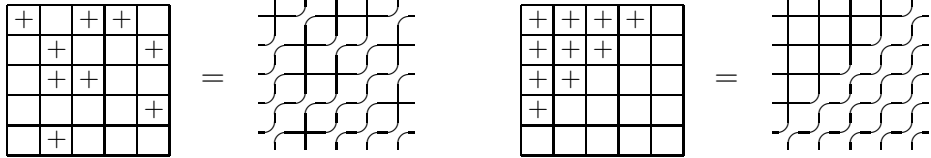
\begin{tinygraph}
\begin{dream}{5}
\hline
   \*  &       &  \*   &  \*   &       \\\hline
       &  \*   &       &       &  \*   \\\hline
       &  \*   &  \*   &       &       \\\hline
       &       &       &       &  \*   \\\hline
       &  \*   &       &       &       \\\hline
\end{dream}
\normal{\quad=\quad}
\begin{pipedream}{5}
\phperm&\phperm&\phperm&\phperm&\phperm\\[-4.5ex]
   \+  &  \jr  &  \+   &  \+   &  \jr  \\
  \jr  &  \+   &  \jr  &  \jr  &  \+   \\
  \jr  &  \+   &  \+   &  \jr  &  \jr  \\
  \jr  &  \jr  &  \jr  &  \jr  &  \+   \\
  \jr  &  \+   &  \jr  &  \jr  &  \jr  \\
\end{pipedream}
\ \qquad\qquad
\begin{dream}{5}
\hline
  \*   &  \*   &  \*   &  \*   &       \\\hline
  \*   &  \*   &  \*   &       &       \\\hline
  \*   &  \*   &       &       &       \\\hline
  \*   &       &       &       &       \\\hline
       &       &       &       &       \\\hline
\end{dream}
\normal{\quad=\quad}
\begin{pipedream}{5}
\phperm&\phperm&\phperm&\phperm&\phperm\\[-4.5ex]
  \+   &  \+   &  \+   &  \+   &  \jr  \\
  \+   &  \+   &  \+   &  \jr  &  \jr  \\
  \+   &  \+   &  \jr  &  \jr  &  \jr  \\
  \+   &  \jr  &  \jr  &  \jr  &  \jr  \\
  \jr  &  \jr  &  \jr  &  \jr  &  \jr  \\
\end{pipedream}
\end{tinygraph}
The square tile boundaries are omitted from the tilings forming the
newtworks of {\em pipes}\/ on right sides of these equalities.  Pipe
dreams are special cases of diagrams introduced by Fomin and Kirillov
\cite{FK96}; for more background, see \cite[Section~1.4]{grobGeom}.

A pipe dream~$P$ yields a word in the Coxeter generators
$s_1,s_2,s_3,\ldots\ $ of~$S_\infty$ by reading the antidiagonal
indices of the crosses in~$P$ along rows, right to left, starting from
the top row and proceeding downward \cite{BB,FK96}.  The {\em Demazure
product}\/ $\delta(P)$ is obtained (as in
\cite[Definition~3.1]{subword}) by omitting adjacent transpositions
that decrease length.  More precisely, $\delta(P)$ is obtained by
multiplying the word of~$P$ using the idempotence relation $s_i^2 =
s_i$ along with the usual braid relations $s_i s_{i+1} s_i = s_{i+1}
s_i s_{i+1}$ and $s_i s_j = s_j s_i$ for $|i-j| \geq 2$.  Up to signs,
this amounts to taking the product of the word of~$P$ in the
degenerate Hecke algebra \cite{FK96}.  Let
\begin{eqnarray*}
  \PP(w) &=& \{\hbox{pipe dreams } P \mid \delta(P) = \wt w\}
\end{eqnarray*}
for a $\kl$ partial permutation~$w$ be the set of pipe dreams whose
Demazure product is the minimal completion of~$w$ to a permutation
\mbox{$\wt w \in S_\infty$}.  Every pipe dream in $\PP(w)$ fits inside
the $\kl$ rectangle, and is to be considered as a pipe dream of size
$\kl$.  The subset of $\PP(w)$ consisting of {\em reduced}\/ pipe
dreams (or {\em rc-graphs}\/ \cite{BB}), where no pair of pipes
crosses more than once, is denoted by~$\rp(w)$.

Here is the observation that will make the limiting arguments in
\cite[Section~6]{quivers} work on nonreduced pipe dreams (see
Proposition~\ref{p:ell}, below).

\begin{lemma} \label{l:dem}
Suppose that~$P \in \PP(w)$.  Then the crossing tiles in~$P$ lie in
the union of all reduced pipe dreams for~$w$.
\end{lemma}
\begin{proof}
The statement is obvious if~$P$ is reduced, so suppose otherwise.
Then some pipe dream $P' \in \PP(w)$ be can be obtained by deleting a
single crossing tile from~$P$.  By induction, every crossing tile
in~$P'$ lies in some reduced pipe dream for~$w$.  On the other hand,
\cite[Theorem~3.7]{subword} implies that a second pipe dream $P'' \in
\PP(w)$ can be obtained from~$P$ by deleting a different crossing
tile.  Induction shows that every crossing tile in~$P''$, including
the tile $P \minus P'$, lies in a reduced pipe dream for~$w$.%
\end{proof}

Lemma~\ref{l:dem} implies that \cite[Corollary~6.10]{quivers} holds as
well for every pipe dream with Demazure product~$v(\rr)$.  This claim
will be made precise in Proposition~\ref{p:ell},~below.

The {\em exponential reverse monomial}\/ associated to a $d \times d$
pipe dream~$P$ is
\begin{eqnarray*}
  (\1-\tilde\xx/\tilde\yy)^P &=& \prod_{+\,\in\,P} (1 - \tilde
  x_+/\tilde y_+),
\end{eqnarray*}
where the variable~$\tilde x_+$ sits at the left end of the row
containing~$\textcross$ after reversing each of the $\xx$~alphabets
before Definition~\ref{d:ratio}, and the variable~$\tilde y_+$ sits
atop the column containing~$\textcross$ after reversing each of the
$\yy$~alphabets
there.  (The row and column labeling in \cite[Section~2.2]{quivers} is
the one meant on the unreversed alphabets here.)

As in \cite[Definition~1.10]{quivers}, let $D_\homv$ be the Ferrers
shape of all locations strictly above the block superantidiagonal.  To
make the meaning of~$\delta(P)$ clear, it is necessary to consider all
crosses in~$P$, including those in $D_\homv$ (unlike the convention of
\cite{quivers}).  Here is the \K-theoretic analogue of
\cite[Proposition~6.9]{quivers}.

\begin{thm}[Pipe formula] \label{t:FK}
The double quiver \K-polynomial is the alternating~sum
\begin{eqnarray*}
  \KQ_\rr(\xx/\oyy) &=& \sum_{\delta(P) = v(\rr)} (-1)^{|P| -
  \length(v(\rr))}(\1-\tilde\xx/\tilde\yy)^{P \minus D_\homv}
\end{eqnarray*}
of exponential reverse monomials associated to pipe dreams $P \minus
D_\homv$ for $P \in \PP(v(\rr))$.  The exponent on $-1$ is the number
crosses in~$P$ minus the length $\length(v(\rr))\!$ of~$v(\rr)$.
\end{thm}
\begin{proof}
Use Definition~\ref{d:ratio} and the symmetry of the double
Grothendieck polynomial for~$v(\rr)$ in each of its $2n+2$ alphabets,
along with the formula of Fomin and Kirillov \cite[Theorem~2.3
and~p.~190]{FK94} (or see \cite[Theorem~4.1 and
Corollary~5.4]{subword}).%
\end{proof}

As in \cite[Section~4.4]{quivers}, let $m+\rr$ be the rank array
obtained from~$\rr$ by adding the nonnegative integer~$m$ to each
entry of~$\rr$.  Let $\xx_{m+\rr}$ be a list of finite alphabets of
sizes $m+r_0,\ldots,m+r_n$, and let the alphabets in
$\oyy{}_{\!m+\rr}$ have sizes \mbox{$m+r_n,\ldots,m+r_0$}.  Denote by
$D_\homv(m)$ the unique reduced pipe dream for the Zelevinsky
permutation $v(m+\rr(\homv))$ in~$S_{d+m(n+1)}$ associated to the
maximal irreducible rank array.

\begin{prop} \label{p:ell}
There is a fixed integer~$\ell$, independent of~$m$, such that for
every pipe dream $P \in \PP(v(m+\rr))$ with at least one
cross~$\textcross$ in an antidiagonal block, setting the last~$\ell$
variables to~$1$ in every finite alphabet from the lists $\xx_{m+\rr}$
and\/~$\oyy{}_{\!m+\rr}$ kills the exponential reverse monomial
$(\1-\tilde\xx/\tilde\yy)^{P\minus D_\homv(m)}$.%
\end{prop}
\begin{proof}
This follows immediately from \cite[Proposition~6.10]{quivers} and
Lemma~\ref{l:dem}.%
\end{proof}

Observe that any pipe dream $P \in \PP(v(\rr))$ with no crossing tiles
in its antidiagonal blocks has its ``interesting'' crosses confined to
the block superantidiagonal.  All other blocks above the antidiagonal
are filled completely with crossing tiles (in \cite{quivers} these are
the $*$ entries), while blocks below the block antidiagonal are empty.
These kinds of pipe dreams $P \in \PP(v(\rr))$ are central to the next
section.

\section{Nonminimal lacing diagrams}

\label{sec:lacing}

Suppose that $\ww = (w_1,\ldots,w_n)$ is a list of partial
permutations in which $w_j$ has size $r_{j-1} \times r_j$.  The
list~$\ww$ can be identified with the (nonembedded) graph in the plane
called its {\em lacing diagram} in \cite[Section~3.1]{quivers}, based
on diagrams of Abeasis and Del\thinspace{}\thinspace{}Fra \cite{AD}.
The vertex set of the graph consists of $r_j$ bottom-justified dots in
column~$j$ for $j = 0,\ldots,n$, with an edge connecting the dot at
height~$\alpha$ (from the bottom) in column~$j-1$ with the dot at
height~$\beta$ in column~$j$ if and only if the entry of~$w_j$ at
$(\alpha,\beta)$ is~$1$.  A~{\em lace}\/ is a connected component of a
lacing diagram.  For example, here is the lacing diagram associated to
a partial permutation list:
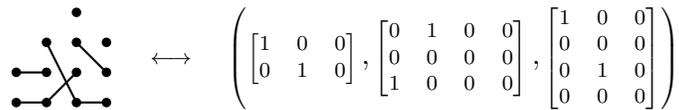
\vspace{-1ex}
\begin{tinygraph}
\pspicture[.1](0,0)(3,3)
\psdots(0,-1)(0,0)(1,-1)(1,0)(1,1)(2,-1)(2,0)(2,1)(2,2)(3,-1)(3,0)(3,1)
\psline(0,-1)(1,-1)(2,0)
\psline(0,0)(1,0)
\psline(1,1)(2,-1)(3,-1)
\psline(2,1)(3,0)
\endpspicture\
\ \quad\longleftrightarrow\quad\
\left(\begin{bmatrix}
  1&0&0 \\
  0&1&0
\end{bmatrix}\normal{,}
\begin{bmatrix}
0&1&0&0\\
0&0&0&0\\
1&0&0&0
\end{bmatrix}\normal{,}
\begin{bmatrix}
1&0&0\\
0&0&0\\
0&1&0\\
0&0&0
\end{bmatrix}\right)
\end{tinygraph}

The goal of this section is to define what it means for a rank array
to equal the Demazure product~$\delta(\ww)$ of a lacing diagram~$\ww$.
That $\delta(\ww)$ is a rank array rather than a minimal lacing
diagram is in analogy with Demazure products of lists of simple
reflections, which are permutations rather than reduced
decompositions.  Usually $\delta(\ww)$ will not equal the rank array
of~$\ww$ itself.  In analogy with Demazure products of {\em reduced}\/
words, however, the Demazure product of a {\em minimal}\/ lacing
diagram will equal its own rank array.

Given a pipe dream~$P$, as in \cite[Theorem~4.4]{subword} say that $P$
{\em simplifies}\/ to \mbox{$D \subseteq P$} if $D$ is the
lexico\-graphically first subword of~$P$ with Demazure product
$\delta(P)$.  Equivalently, denoting by $P_{\leq i}$ the length~$i$
initial string of simple reflections in~$P$, the simplification $D$ is
obtained from~$P$ by omitting the~$i^\th$ reflection from~$P$ for
all~$i$ such that $\delta(P_{\leq i-1}) = \delta(P_{\leq i})$.

\begin{lemma} \label{l:simp}
Suppose that $P$ is a $\kl$ pipe dream and let $\alpha|P$ be
the $k \times (\ell + \alpha)$ pipe dream obtained by adding $\alpha$
columns of~$k$~$\textcross$ tiles to the left side of~$P$.  The pipe
dream $P$ simplifies to~$D$ if and only if $\alpha|P$ simplifies
to~$\alpha|D$.
\end{lemma}
\begin{proof}
Since $P$ is reduced if and only if $\alpha|P$ is reduced, we may
assume that $P$ is not reduced.  Moreover, by adding $\textcross$
tiles to~$P$ one by one (from right to left in each row and top to
bottom, as usual), it is enough to prove the lemma when $|P| = 1 +
|D|$.  In this case, a single pair of pipes in~$P$ crosses twice, as
does the corresponding pair of pipes (shifted to the right
by~$\alpha$) in~$\alpha|P$.  The simplifications of~$P$ and $\alpha|P$
are obtained by deleting the southwestern crossings of the
corresponding pairs of pipes.%
\end{proof}

\begin{defn}
Suppose $P_1,\ldots,P_n$ are pipe dreams of sizes $r_0 \times r_1,
\ldots, r_{n-1} \times r_n$, and set $d = r_0 + \cdots + r_n$.  Denote
by $P(P_1,\ldots,P_n)$ the $d \times d$ pipe dream in which every
block strictly above the block superantidiagonal is filled with
crossing tiles, and the superantidiagonal $r_{j-1} \times r_j$ block
in block row~\mbox{$j-1$} is the pipe dream~$P_j$.
\end{defn}

Given a $\kl$ pipe dream~$P$, let $\oP$ be the $\kl$ pipe dream that
results after rotating~$P$ through~$180^\circ$.  Also, recall from
\cite[Theorem~3.7]{BB} the notion of {\em top}\/ pipe dream for a
partial permutation~$w$, which is the unique reduced pipe dream
in~$\rp(w)$ that has no elbow tile due north of a crossing tile.

\begin{prop} \label{p:dem}
Fix a lacing diagram $\ww = (w_1,\ldots,w_n)$.  The Demazure product
of $P(P_1,\ldots,P_n)$ is independent of $P_1,\ldots,P_n$, as long as
$\oP_j \in \PP(w_j)$ for all $j = 1,\ldots,n$.
\end{prop}
\begin{proof}
Associativity of Demazure products implies that we can take Demazure
products first in each block row.  By Lemma~\ref{l:simp} these
stripwise Demazure products don't change when each $P_j$ is replaced
by its simplification.  Neither do the Demazure products of the pipe
dreams~$\oP_j$.  Therefore we can assume that each $P_j$---and hence
each block row of $P(P_1,\ldots,P_n)$---is reduced.

The Demazure product of each block row is unchanged by chute and
inverse chute moves \cite{BB} that remain within block rows, because
the Demazure product equals the usual product on reduced expressions.
In addition, under these operations the reduced pipe dreams~$\oP_j$
remain inside of $\rp(w_j)$ for all~$j$.  Therefore
$\delta(P(P_1,\ldots,P_n))$ equals the Demazure product of the pipe
dream $P(D_1,\ldots,D_n)$ in which $\check{D}_j$ is the unique ``top''
reduced pipe dream for~$w_j$ by \cite[Theorem~3.7]{BB}.%
\end{proof}

\begin{defn}
Fix a lacing diagram~$\ww$.  If, for some (and hence, by
Proposition~\ref{p:dem}, every) sequence $P_1,\ldots,P_n$ of pipe
dreams satisfying $\oP_j \in \PP(w_j)$ for all~$j$, the Demazure
product of $P(P_1,\ldots,P_n)$ is a Zelevinsky permutation~$v(\rr)$,
then we write $\delta(\ww) = \rr$ and call the rank array~$\rr$ the
\bem{Demazure product} of the lacing diagram~$\ww$.
\end{defn}

\section{Rank stability of lacing diagrams}

\label{sec:stability}

Next we show that lacing diagrams with Demazure product~$\rr$ are
stable, in the appropriate sense, under uniformly increasing ranks
obtained by replacing $\rr$ with~\mbox{$m+\rr$}.  To ease the
language, we use `horizontal strip~$j$' as a synonym for
`block~row~$j$'.

\begin{lemma} \label{l:fits}
If $P(P_1,\ldots,P_n) \in \PP(v(1+\rr))$ and each $\oP_{\!j}$ is the
top pipe dream for a $(1+r_{j-1}) \times (1+r_j)$ partial
permutation~$w_j$, then all crossing tiles of $P_j$ lie in the
southwest $r_{j-1} \times r_j$ rectangle of the antidiagonal block
in horizontal strip~\mbox{$j-1$}.
\end{lemma}

Thus the antidiagonal block in the Lemma is supposed to have one blank
row on top and one blank column to the right of the southwest $r_{j-1}
\times r_j$ rectangle in question.

\begin{proof}
No reduced pipe dream for~$v(1+\rr)$ has a crossing tile on the main
superantidiagonal, by \cite[Proposition~5.15]{quivers}.
Lemma~\ref{l:dem} implies that the same is true of~$P$.  It follows
that $w_j = 1 + w_j'$ for some $r_{j-1} \times r_j$ partial
permutation~$w_j'$.  Consequently, the left column of~$\oP_j$ has no
crossing tiles, and shifting all crossing tiles in~$\oP_j$ one unit to
the left results in the top pipe dream for~$w_j'$.  This top pipe
dream fits inside the rectangle of size $r_{j-1} \times r_j$.%
\end{proof}

Suppose $P = P(P_1,\ldots,P_n)$ is a pipe dream in which
\begin{equation} \tag{SW}
\begin{array}{@{\ }l@{}}
\hbox{$P_j$ has size $(1+r_{j-1}) \times (1+r_j)$, but
	every~$\textcross$ in~$P_j$ lies in the}
\\
\hbox{southwest $r_{j-1} \times r_j$ rectangle.}
\end{array}
\end{equation}%
Write $P_j'$ for the $r_{j-1} \times r_j$ pipe dream consisting of the
southwest rectangle of~$P_j$, and then write $P' =
P(P_1',\ldots,P_n')$.  Thus $P$ has block sizes consistent with
ranks~$1+\rr$, while $P'$ has block sizes consistent with ranks~$\rr$.
The construction can also be reversed to create~$P$ having been given
the pipe dream called~$P'$.

Given a reduced pipe dream~$D$, an elbow tile is {\em absorbable}\/
\cite[Section~4]{subword} if the two pipes passing through it
intersect in a crossing tile to its northeast.  It follows from the
definitions that a pipe dream~$P$ simplifies to~$D$ if and only if $P$
is obtained from~$D$ by changing (at will) some of its absorbable
elbow tiles into crossing tiles.

\begin{lemma} \label{l:abs}
Suppose $D = (D_1,\ldots,D_n)$ satisfies the (SW) condition.  Then $D$
is a reduced pipe dream for~$v(1+\rr)$ if and only if $D' =
(D_1',\ldots,D_n')$ is a reduced pipe dream for~$v(\rr)$.  In this
case, the absorbable elbow tiles in horizontal strip~\mbox{$j-1$}
of~$D'$ are in bijection with the absorbable elbow tiles in the
southwest $r_{j-1} \times r_j$ rectangle of the antidiagonal block in
\mbox{horizontal strip~$j-1$ of~$D$}.
\end{lemma}
\begin{proof}
The first claim is a straightforward consequence of
\cite[Proposition~5.15]{quivers}.  The second claim follows because
the corresponding pairs of pipes in~$D$ and~$D'$ pass through
corresponding elbow tiles.  The rest of the proof makes this statement
precise.

Given a nonzero entry of the Zelevinsky permutation~$v(1+\rr)$,
exactly one of the following three conditions must hold: (i)~the entry
lies in the northwest corner of some superantidiagonal block; (ii)~the
entry lies in the southeast corner of the whole matrix; or (iii)~there
is a corresponding nonzero entry in~$v(\rr)$.  This means that the
pipes in~$D'$ are in bijection with those pipes in~$D$ corresponding
to nonzero entries of~$v(1+\rr)$ that do not satisfy (i) or~(ii).
Furthermore, it is easily checked that the pipes in~$D$ of type~(i)
or~(ii) can only intersect a superantidiagonal block in its top row or
rightmost column.  Hence to say
\begin{quote}
the two pipes passing through an elbow tile in the southwest
\mbox{$r_{j-1} \times r_j$} rectangle of the antidiagonal block in
horizontal strip~$j-1$ of~$D$ correspond to the pipes passing through
the corresponding elbow tile in~$D'$
\end{quote}
actually makes sense.  That this claim is true follows from
\cite[Proposition~5.15]{quivers}, and it immediately proves the
lemma.%
\end{proof}

Here is the \K-theoretic (nonminimal lacing diagram) analogue of
\cite[Corollary~5.16]{quivers}.  The notation is as in
\cite[Section~4.4]{quivers}: given $m \in \NN$ and a partial
permutation~$w$, the partial permutation $m+w$ is obtained by
letting~$w$ act on $m + \ZZ_{> 0} = \{m+1,m+2,\ldots\}$ in the obvious
manner instead of on $\ZZ_{> 0} = \{1,2,\ldots\}$.  For a list $\ww =
(w_1,\ldots,w_n)$ of partial permutations, set $m + \ww =
(m+w_1,\ldots,m+w_n)$.

\begin{prop} \label{p:L}
For each array~$\rr$, let $L(\rr) = \{\ww \mid \delta(\ww) = \rr\}$ be
the set of lacing diagrams~$\ww$ with Demazure product\/~$\rr$.  Then
$L(\rr)$ and $L(m+\rr)$ are in \mbox{canonical bijection}:
\begin{eqnarray*}
  L(m+\rr) &=& \{m+\ww \mid \ww \in L(\rr)\}.
\end{eqnarray*}
\end{prop}
\begin{proof}
It suffices to prove the case $m = 1$, so suppose $\ww \in L(1+\rr)$.
Let $P = P(P_1,\ldots,P_n)$ be the pipe dream in $\PP(v(1+\rr))$ for
which each~$\oP_j$ is the top pipe dream in~$\rp(w_j)$.  Then $P$
simplifies to a reduced pipe dream $D \in \rp(v(1+\rr))$.  By
Lemma~\ref{l:fits} there is a corresponding pipe dream $D' \in
\rp(v(\rr))$, constructed via the procedure after Lemma~\ref{l:fits}.
On the other hand, the pipe dream~$P'$ constructed from~$P$ results by
changing back into crossing tiles those elbow tiles in~$D'$ that
correspond to the $\textcross$ tiles deleted from~$P$ to get~$D$.
Lemma~\ref{l:abs} says that $P'$ has Demazure product~$v(\rr)$.
Defining $\ww'$ by the equality $1+\ww' = \ww$, which can be done by
Lemma~\ref{l:fits}, it follows that $\ww' \in L(\rr)$.

In summary, we have constructed $P'$ from~$P$ via the intermediate
steps
$$
  P \in \PP(v(1+\rr))\ \goesto\ D \in \rp(v(1+\rr))\ \goesto\ D' \in
  \rp(v(\rr))\ \goesto\ P' \in \PP(v(\rr)),
$$
where the first and third steps are simplification and
``unsimplification''.  Consequently, $L(1+\rr) \subseteq \{1+\ww' \mid
\ww' \in L(\rr)\}$.  But the arguments justifying these steps are all
reversible, so the reverse containment holds, as well.%
\end{proof}

\section{Stable double component formula}

\label{sec:formula}

The main result in this paper, namely Theorem~\ref{t:formula},
involves {\em stable double Grothen\-dieck polynomials}\/
$\hat\GG_w(\zz/\dot\zz)$ for $\kl$ partial permutations~$w$
\cite{FK94}, which we recall presently.  Suppose that the argument of
a Laurent polynomial~$\GG$ is naturally a pair of alphabets
$\zz$~and~$\dot\zz$ of sizes $k$ and~$\ell$, respectively.  In this
section and the next, the convention is that if $\GG(\zz/\dot\zz)$ is
written, but $\zz$ or $\dot\zz$ has {\em fewer}\/ than the required
number of letters, then the rest of the letters are assumed to
equal~$1$.  For example, the notation
$\KQ_{m+\rr}(\xx_\rr/\oyy{}_{\!\rr})$ indicates that all variables in
$\xx_{m+\rr} \minus \xx_\rr$ and $\oyy{}_{\!m+\rr} \minus
\oyy{}_{\!\rr}$ (see the paragraph preceding Proposition~\ref{p:ell})
are to be set equal to~$1$.

Under this convention, let $w$ be a $\kl$ partial permutation, and
write $\GG_{m+w}(\zz_k/\dot\zz_\ell)$ for each $m \geq 0$ to mean the
Laurent polynomial $\GG_{m+w}$ applied to alphabets $\zz_k$
and~$\dot\zz_\ell$ of fixed sizes $k$ and~$\ell$.  As~$m$ gets large,
these Laurent polynomials eventually stabilize, allowing the notation
$\hat\GG_w(\zz/\dot\zz) = \lim_{m \to \infty}
\GG_{m+w}(\zz_k/\dot\zz_\ell)$ for the stable double Grothendieck
polynomial.

Given a lacing diagram~$\ww$ with $r_j$ dots in column~$j$, for $j =
0,\ldots,n$ denote by
\begin{eqnarray*}
  \GG_\ww(\xx/\oyy) &=& \GG_{w_1}(\xx^0/\yy^1)\cdots
  \GG_{w_n}(\xx^{n-1}/\yy^n)
\end{eqnarray*}
the product of double Grothendieck polynomials taken over partial
permutations in the list $\ww = (w_1,\ldots,w_n)$.  Add hats over
every~$\GG$ for the stable Grothendieck case.

Here now is the main result, the \K-theoretic analogue of the
(cohomological) component formula for stable double quiver polynomials
\cite[Theorem~6.20]{quivers}.

\begin{thm} \label{t:formula}
The limit of double quiver \K-polynomials
$\KQ_{m+\rr}(\xx_\rr/\oyy{}_{\!\rr})$ for\/ $m$ approaching\/~$\infty$
exists and equals the alternating sum
\begin{eqnarray*}
  \fG_\rr(\xx/\oyy) \ \::=\ \:
  \lim_{m\to\infty}\KQ_{m+\rr}(\xx_\rr/\oyy{}_{\!\rr}) &=&
  \!\!\!\sum_{\ww \in L(\rr)}
  (-1)^{\length(\ww)-d(\rr)}\hat\GG_\ww(\xx/\oyy)
\end{eqnarray*}
of products of stable double Grothendieck polynomials, where $L(\rr) =
\{\ww \mid \delta(\ww) = \rr\}$, $\length(\ww) = \sum_{i=1}^n
\length(\wt w_i)$, and $d(\rr) = \length(v(\rr)) - \length(v(\homv))$.
The limit polynomial~$\fG_\rr(\xx/\oyy)$ is symmetric separately in
each of the $2n+2$ finite alphabets
$\xx^0,\ldots,\xx^n,\yy^n,\ldots,\yy^0$.%
\end{thm}

\begin{defn} \label{d:stableK}
$\fG_\rr(\xx/\oyy)$ is called the \bem{stable double quiver
\K-polynomial}.%
\end{defn}

As we shall see in Corollary~\ref{c:formula} and the comments after
it, the Laurent polynomial $\fG_\rr(\xx/\oyy)$ is not a new object: it
is obtained from Buch's power series~$P_r$ \cite[Section~4]{Buch02} by
substituting $1-x_i$ for~$x_i$ and $1-y_j^{-1}$ for~$y_j$ in each
polynomial~$G_{\mu_k}$~there.

\begin{proof}
Define $\KQ_{m+\rr}(\xx/\oyy)_\ell$ by setting the last~$\ell$
variables to~$1$ in every finite alphabet from the lists $\xx_{m+\rr}$
and\/~$\oyy{}_{\!m+\rr}$.  Similarly, for each lacing diagram~$\ww$,
define $\GG_{m+\ww}(\xx/\oyy)_\ell$ by setting the same variables
to~$1$ in $\GG_{m+\ww}(\xx/\oyy)$.  Because of the nature of the limit
in question, and the defining properties of stable Grothendieck
polynomials, it suffices to prove that for all $m \geq 0$ and some
fixed~$\ell$ independent of $m$,
\begin{eqnarray*}
  \KQ_{m+\rr}(\xx/\oyy)_\ell &=& \!\!\!\sum_{\ww \in L(\rr)}
  (-1)^{\length(\ww)-d(\rr)}\GG_{m+\ww}(\xx/\oyy)_\ell.
\end{eqnarray*}

Fix $\ell$ as in Proposition~\ref{p:ell}, and apply Theorem~\ref{t:FK}
to $m+\rr$ instead of~$\rr$.  Setting the last $\ell$ variables in
each alphabet to~$1$ on the right hand side there kills all summands
corresponding to pipe dreams~$P$ that are not expressible as
$P(P_1,\ldots,P_n)$ for some list of pipe dreams~$P_j$ of sizes
$(m+r_{j-1}) \times (m+r_j)$; this is the content of
Proposition~\ref{p:ell}.  What remains on the right side of
Theorem~\ref{t:FK} is a sum of terms having the form $(-1)^{|P| -
\length(v(m+\rr))}(1-\tilde\xx/\tilde\yy)^{P \minus D_\homv(m)}_\ell$
for pipe dreams~$P = P(P_1,\ldots,P_n)$ in $\PP(v(m+\rr))$.  If $P_j
\in \PP(m+w_j)$ for each~$j$, then this term equals the product
\begin{stareqn} \label{*}
   &\displaystyle(-1)^{\length(\ww)-d(\rr)}\prod_{i=1}^n
  (-1)^{|P_j| - \length(\wt w_j)}
  (1-\tilde\xx^{j-1}/\tilde\yy^j)^{P_j}_\ell&
\end{stareqn}%
for $\ww = (w_1,\ldots,w_n)$.  The signs in~\eqref{*} are correct
because $|P| - \length(v(m+\rr)) = \sum_j |P_j| - d(m+\rr)$, and
$d(m+\rr) = d(\rr)$.  To make sense of
$(1-\tilde\xx^{j-1}/\tilde\yy^j)^{P_j}$, identify $P_j$ with the $d
\times d$ pipe dream consisting of just $P_j$ on the $j^\th$
superantidiagonal block.

For each lacing diagram $\ww \in L(\rr)$, let $\PP_\ww(m+\rr)$ be the
set of pipe dreams $P(P_1,\ldots,P_n) \in \PP(v(m+\rr))$ such that
$\oP_j \in \PP(m+w_j)$ for all~$j$.  Summing the products in~\eqref{*}
over pipe dreams $P \in \PP_\ww(m+\rr)$ yields
$(-1)^{\length(\ww)-d(\rr)}\GG_{m+\ww}(\xx/\oyy)_\ell$ by
\cite[Theorem~2.3 and~p.~190]{FK94} (see also
\cite[Section~5]{subword}).  Summing over $\ww \in L(\rr)$ completes
the proof, by Proposition~\ref{p:L}.%
\end{proof}

\begin{remark}
Theorem~\ref{t:formula} implies that
$\fG_{m+\rr}(\xx_\rr/\oyy{}_{\!\rr}) = \fG_\rr(\xx/\oyy)$, in analogy
with the (defining) stability properties of stable double Grothendieck
polynomials.
\end{remark}

\begin{remark}
Theorem~\ref{t:formula} gives an explicit combinatorial formula, but
the characterization of the Demazure product $\delta(\ww)$ of a lacing
diagram via Zelevinsky permutations would be more satisfying if it
were intrinsic.  That is, it would be better to identify those partial
permutation lists that fit stripwise into a pipe dream with Demazure
product~$v(\rr)$ using the language of lacing diagrams, without
referring to Zelevinsky permutations or pipe dreams.  Such an
intrinsic method appears in \cite{BFR03}.
\end{remark}

\section{Sign alternation}

\label{sec:alt}

A permutation $\mu \in S_\infty$ is {\em grassmannian}\/ if it has at
most one descent---that is, at most one index~$p$ such that $\mu(p) >
\mu(p+1)$.  A~crucial property of arbitrary stable double Grothendieck
polynomials, proved in \cite[Theorem~6.13]{BuchLR}, is that every such
polynomial $\hat\GG_w(\zz/\dot\zz)$ has a unique expression
\begin{eqnarray*}
  \hat\GG_w(\zz/\dot\zz) &=& \sum_{{\rm grassmannian}\ \mu}
  \alpha_w^\mu\hat\GG_\mu(\zz/\dot\zz)
\end{eqnarray*}
as a sum of stable Grothendieck polynomials~$\hat\GG_\mu$ for
grassmannian permutations.  If $\ub\mu = (\mu_1,\ldots,\mu_n)$ is a
sequence of partial permutations such that the minimal completions
$\wt \mu_1,\ldots,\wt \mu_n$ are grassmannian, then let us
call~$\ub\mu$ a {\em grassmannian}\/ lacing diagram.

\begin{cor} \label{c:formula}
If $\alpha_\ww^\ub\mu = \prod_{i=1}^n \alpha_{w_i}^{\mu_i}$ for each
lacing diagram~$\ww$ and grassmannian~$\ub\mu$,~then
\begin{eqnarray*}
  \fG_\rr(\xx/\oyy) &=&\ \sum_\ub\mu c_\ub\mu(\rr)
  \hat\GG_\ub\mu(\xx/\oyy)\\
\makebox[0pt][r]{for the constants\qquad}
  c_\ub\mu(\rr) &=& \sum_{\ww \in L(\rr)}
  (-1)^{\length(\ww)-d(\rr)}\alpha_\ww^\ub\mu,
\end{eqnarray*}
where the first sum above is over all grassmannian lacing
diagrams~$\ub\mu$.
\end{cor}
\begin{proof}
Expand the right hand side of Theorem~\ref{t:formula} using $\hat\GG_w
= \sum_\mu\alpha_w^\mu\hat\GG_\mu$.%
\end{proof}

Let $\fG_\rr(\xx/\oxx)$ be the specialization of the stable double
quiver \K-polynomial obtained by setting $\yy^j = \xx^j$ for $j =
0,\ldots,n$.  Independently from Corollary~\ref{c:formula}, it follows
from \cite[Theorem~4.1]{Buch02} that the (ordinary) stable quiver
\K-polynomial
\begin{eqnarray*}
  \fG_\rr(\xx/\oxx) &=& \sum_{\ub\mu} c_\ub\mu(\rr)
  \hat\GG_\ub\mu(\xx/\oxx)
\end{eqnarray*}
is a sum of products of stable double Grothendieck polynomials
$\hat\GG_{\mu_j}(\xx^{j-1}/\xx^j)$ for grassmannian permutations~$\wt
\mu_j$, with uniquely determined integer coefficients~$c_\ub\mu(\rr)$.
That these coefficients are the same as in Corollary~\ref{c:formula}
follows from the fact that the right side
above determines the same element in the $n^\th$ tensor power of
Buch's bialgebra~$\Gamma$ from \cite{BuchLR,Buch02} as does the right
side of the top formula in Corollary~\ref{c:formula}.

In addition to proving
the expansion of~$\hat\GG_w$ as a sum of terms
$\alpha_w^\mu\hat\GG_\mu$, Buch showed in \cite[Theorem~6.13]{BuchLR}
that the coefficients $\alpha_w^\mu$ can only be nonzero if
$\length(\mu) \geq \length(w)$, and he conjectured that the sign
of~$\alpha_w^\mu$ equals $(-1)^{\length(\mu)-\length(w)}$.  This was
proved by Lascoux \cite[Theorem~4]{Las01} as part of his extension of
``transition'' from Schubert polynomials to Grothendieck polynomials.
Since, as shown in \cite[Section~5]{Buch02}, the coefficients
$\alpha_w^\mu$
are special cases of the coefficients~$c_\ub\mu(\rr)$, Lascoux's
result is evidence for the following more general statement that was
surmised by Buch (prior to \cite{Las01}).

\begin{thm}[{\cite[Conjecture~4.2]{Buch02}}]\label{t:alt}
The coefficients $c_\ub\mu(\rr)$ alternate in sign; that is,
\mbox{$(-1)^{\length(\ub\mu)-d(\rr)} c_\ub\mu(\rr) \geq 0$} is a
nonnegative integer.
\end{thm}
\begin{proof}
By \cite[Theorem~4]{Las01} the sign of~$\alpha_\ww^\ub\mu$ is
$(-1)^{\length(\ub\mu)-\length(\ww)}$.  Thus the sign
of~$c_\ub\mu(\rr)$~is $(-1)^{\length(\ww) -
d(\rr)}(-1)^{\length(\ub\mu) - \length(\ww)} = (-1)^{\length(\ub\mu) -
d(\rr)}$, by the second formula in Corollary~\ref{c:formula}.%
\end{proof}


\def\cprime{$'$}
\providecommand{\bysame}{\leavevmode\hbox to3em{\hrulefill}\thinspace}


\end{document}